\documentclass{article}
\usepackage{spconf,amsmath,graphicx}
\usepackage{amssymb,amsfonts}
\usepackage{algorithm,algorithmic}

\title{ New Accumulative Score Function based Bound for Sparsity Level of \quad $l_1$ Minimization}
  \name{$Sheng\ Han^{1*}, Suzhen\ Wang^2\sthanks{S. Han and S.Z. Wang contributed equally to this work}, Zhiguo\ Zhang^1$}
	\address {$^1$Department of Electrical and Electronic Engineering\\
	 The University of Hong Kong, Hong Kong\\
   $^2$Department of Information Engineering\\
	The Chinese University of Hong Kong, Hong Kong\\
     }
\begin{document}
%
\maketitle
\begin{abstract}
  This paper discusses a fundamental problem in compressed sensing: the sparse recoverability of $l_1$ minimization with an arbitrary sensing matrix. We develop an new accumulative score function (ASF) to provide a lower bound for the recoverable sparsity level ($SL$) of a sensing matrix while preserving a low computational complexity. We first define a score function for each row of a matrix, and then ASF sums up large scores until the total score reaches 0.5. Interestingly, the number of involved rows in the summation is a reliable lower bound of $SL$. It is further proved that ASF provides a sharper bound for $SL$ than coherence We also investigate the underlying relationship between the new ASF and the classical RIC and achieve a RIC-based bound for $SL$.
\end{abstract}
\section{Introduction}
\label{sec:intro}
Recently, compressed sensing (CS) has become a powerful technique for exploring sparse representation of a signal given a redundant dictionary. In CS, it is a fundamental problem to study the sparsity level $(SL)$ of $l_1$ minimization. Specifically speaking, $ s\le SL$ means for any $s$-sparse vector $x$ (i.e., $\|x\|_0\le s$) can be correctly recovered by solving the following $l_1$ minimization problem.
\begin{equation}
\nonumber
(L_1)\quad \quad  \hat{x} = \arg\min_{x\in R^{m}}\|x\|_1, ~ s.t. \ y = Ax,
\end{equation}
where $A\in R^{n\times m}$ ($n<m$) is a sensing matrix, $y\in R^{n}$ is the sensed data and $\hat{x}\in R^m$ is the $l_1$ minimizer.
Generally, the larger the $l_0$ norm of the vector $x$, the less possible to precisely recover this vector by solving a corresponding $l_1$ minimization problem. Therefore, we are more interested in providing a precise lower bound of $SL$. Researchers have explored many tools to study the sparse recovery property of the $l_1$ minimization problem. Some useful and powerful tools include null space property (NSP) \cite{gribonval2003sparse,zhang2013theory}, coherence \cite{gribonval2003sparse,donoho2003optimally}, restricted isometry property (RIP) \cite{candes2005decoding,candes2008restricted}, which will be briefly introduced below.

NSP states that $l_1$-minimization holds a sparsity level of $s$ if and only if the following condition holds
\begin{equation}
\nonumber
(NSP)\quad \sum_{i\in \mathcal{S}} |v_i| < \frac{1}{2}\|v\|_1,  \ \forall v\neq 0 \in ker A,
\end{equation}
where $|\mathcal{S}| =s$ and $ker A$ means kernel of matrix $A$. NSP itself is NP-hard but it reveals the fact that the sparse recovery ability of the $l_1$ minimization problem is actually determined by its sensing matrix. NSP has also been widely applied to derive other useful sparsity-related properties on a sensing matrix, such as coherence and RIP.

Coherence is an important measure to get a lower bound of $SL$ as shown below \cite{gribonval2003sparse, donoho2003optimally}:
\begin{equation}
\nonumber
\noindent \quad SL\ge \frac{1}{2}(1+\mu^{-1}) -1,
\end{equation}
where $\mu$ denotes the coherence of $A$. Usually the smaller the coherence, the higher the $SL$. However, coherence may not be sufficient to estimate a lower bound of $SL$ because, as pointed by \cite{candes2008enhancing}, rescaling the columns of a sensing matrix may enhance the sparsity recovery ability of $l_1$ minimization in some cases but does not change the value of the coherence.

Restricted isometry constant (RIC) is a more sophisticated sparsity-relevant constant than coherence \cite{foucart2013mathematical}. The RIP condition of order $s$ states that, there exists a RIC, $\sigma_s\in (0,1)$, that makes the following inequality holds.
\begin{equation}
\nonumber
(RIP)\quad (1-\sigma_s)\|x\|_2^2 \le \|Ax\|_2^2 \le (1+\sigma_s)\|x\|_2^2,
\end{equation}
for any $s$-sparse vector $x\in R^{m}$. $\sigma_s$ increases with the sparsity number $s$. In order to guarantee the correct recovery of a sparse vector via $l_1$ minimization, RIC is required to be less than a certain bound. Tremendous effort has been made to sharpen this bound \cite{cai2013sharp}. But it is still combinatorial complexity to compute RIC of a certain sensing matrix.

In this paper, we introduce a new method, accumulative score function (ASF), to analyze the sparsity recovery ability of the $l_1$ minimization. With the help of NSP, we derive a lower bound of $SL$ based on ASF. ASF considers not only the angle information between columns but also the scale information of each column in a sensing matrix. Importantly, we prove that ASF tends to provide a sharper bound for $SL$ than coherence, and, thus, generally RIP is a more sophisticated tool to study sparse recovery of the $l_1$ minimization. Furthermore, we investigate the underlying relationship between the new ASF and the classical RIC, and derive a new RIC-based formulation for the lower bound of $SL$.

The rest of this paper is organized as follows. Section 2 introduces the ASF and analyze its performance in sparsity estimation. Section 3 is devoted to prove that ASF-based bound is relatively sharper than the coherence-based bound for $SL$, and to derive an alternative RIC-based bound for $SL$. We finally conclude this paper in Section 4.

\section{Accumulative Score Function}

\label{sec:format}
\noindent \textbf{Notations}
Let $[m]$ denote the set $\{1,...,m\}$ and let $\mathcal{T}$ be a subset of $[m]$. Let $\mathcal{T}^c$ denote the complementary set of $\mathcal{T}^c$.
Use $|\mathcal{T}|$ denote the cardinality of a set $\mathcal{T}$.
For a matrix $B\in R^{n\times m}$, $B_{\mathcal{T}}$ denotes a submatrix constructed by columns of matrix $B$ indexed by elements in $\mathcal{T}$.
$\Lambda_{max}(B)$ and $\Lambda_{min}(B)$ respectively denote the maximal and minimal eigenvalue of $B$.

Let $\{\alpha_1,..., \alpha_m\}$ denote columns of the measurement matrix $A\in R^{n\times m}$.
Let $C=A^TA \in R^{m\times m}$. Obviously, the diagonal entries of matrix $C$, denoted as $c_{ii}$, equal to $\|\alpha_i\|^2_2$ for $i\in [m]$.
The non-diagonal entries, $c_{ij}$  $(i\neq j)$, are equal to inner products between $\alpha_i$ and $\alpha_j$ (i.e. $c_{ij}  = <\alpha_i,\alpha_j>$).

We define the score function $\rho(i)$ for each row $i$ of matrix $C$ as:
\begin{equation}
\begin{split}
\label{Proposed Property2}
\rho(i) = \frac{\nu(i)}{\nu(i)+1}, \\
\end{split}
\end{equation}
where
\begin{equation}
\label{Property1}
\nu(i) = \max_{j} \frac{|c_{ij}|}{c_{ii}}, \ \forall j \neq i.
\end{equation}

By summing up the first $s$ largest scores, we define the ASF as
\begin{equation}
\nonumber
\rho(\mathcal{S^*}) = \max_{|\mathcal{S}| =s}\sum_{i\in \mathcal{S}}\rho(i),
\end{equation}
where $\mathcal{S^*}$ denotes optimal index set that makes $\sum_{i\in \mathcal{S}}\rho(i)$ largest.
Without loss of generality, we can assume sequence $\{\rho(1),...,\rho(m)\}$ is in a non-increasing order, that is
\begin{equation}
\begin{split}
\nonumber
\rho(i)\ge \rho(j), \ \forall \ i\le j.
\end{split}
\end{equation}
Then $\rho(\mathcal{S^*})$ with $|\mathcal{S^*}| =s$ can be re-written as
\begin{equation}
\rho(\mathcal{S^*}) = \sum_{i=1}^s \rho(i)
\end{equation}

\subsection{Sparsity Analysis}
We now aim to show how to apply ASF to obtain a lower bound of sparsity level of the $l_1$ minimization problem with an arbitrary sensing matrix.
First we show a lemma which can be taken as an alternative interpretation of NSP.

\noindent \textbf{Lemma 1} Let $\mathcal{S}$ denote the support set of a sparse vector $x$, i.e. the index set of nonzero-entry of $x$, while $\hat{x}$ denotes an $l_1$ minimizer,
then we have
\begin{equation}
\label{Lemma1}
\|x - \hat{x}\|_1 \le 2 \|(x - \hat{x})_S\|_1.
\end{equation}

\noindent \textbf{Proof} Firstly we have $\|x \|_1 \ge \|\hat{x}\|_1$ which means
\begin{equation}
\nonumber
\|(x- \hat{x})_{\mathcal{S}} \|_1\ge \|x\|_1 -\|\hat{x}_{\mathcal{S}}\|_1 \ge \|\hat{x}\|_1 - \|\hat{x}_{\mathcal{S}}\|_1 = \|\hat{x}_{{\mathcal{S}}^c}\|_1.
\end{equation}
Let $e = x -\hat{x}$, then we have
\begin{equation}
\nonumber
\|e\|_1 = \|e_{\mathcal{S}}\|_1 +\|\hat{x}_{{\mathcal{S}}^c}\|_1 \le \|e_{\mathcal{S}}\|_1+ \|(x - \hat{x})_{\mathcal{S}} \|_1 = 2\|e_{\mathcal{S}}\|_1.
\end{equation}
Thus Lemma 1 is proven. $\blacksquare$

Based on the above lemma, we can then show how to apply ASF to assess the sparsity level of the $l_1$ minimization problem in the following theorem.

\noindent \textbf{Theorem 1} Suppose $\mathcal{S}$ is the support set of a sparse vector $x$. If $\rho(\mathcal{S}) < \frac{1}{2}$, $x$ can be correctly
recovered by solving the $l_1$ minimization problem.

\noindent \textbf{Proof} Assume $x\neq \hat{x}$ and let $e = x - \hat{x}$, then we have $A^TAe = 0$ since $Ax= A\hat{x}$.
For each row $i$, we have
\begin{equation}
\nonumber
\sum_{j\neq i} c_{ij}e_j + c_{ii} e_i =0,
\end{equation}
from which we can derive the following inequality:
\begin{equation}
\nonumber
|c_{ii}e_i|\le \sum_{j\neq i}|c_{ij}e_j|\le |c_{ik}|\sum_{j\neq i} |e_j|,
\end{equation}
where $|c_{ik}| = \max_{j\neq i} \{|c_{ij}|\}$.
Then we get
\begin{equation}
\nonumber
(c_{ii} +|c_{ik}|) |e_i| \le |c_{ik}|\|e\|_1.
\end{equation}
Furthermore, we have
\begin{equation}
\label{ElementaryInequality}
|e_i| \le \frac{|c_{ik}|}{c_{ii} + |c_{ik}|}\|e\|_1 =\rho(i)\|e\|_1
\end{equation}
Combining ($\ref{Lemma1}$) in lemma 1 with the inequality in $(\ref{ElementaryInequality})$, we get
\begin{equation}
\nonumber
\frac{1}{2}\|e\|_1\le \|e_{\mathcal{S}}\|_1 \le \rho(\mathcal{S}) \|e\|_1< \frac{1}{2}\|e\|_1,
\end{equation}
which is contradictory. Then we have $x = \hat{x}$. Thus we prove the theorem. $\blacksquare$

It should be noted that $\rho(\mathcal{S})< \frac{1}{2}$ is a sufficient condition to guarantee that the $l_1$ minimization problem can recover the correct sparse vector. In general, information about the support set of a sparse vector is unknown. But we can still use ASF to provide some useful information about the sparsity level of a sensing matrix. The following corollary presents a simple way to get a lower bound of $SL$ via ASF.

\noindent \textbf{Corollary 1} Suppose sequence $\{\rho(1),...,\rho(m)\}$ is in a non-increasing order, then we have $SL\ge l^*$, where $l^*$ is determined below:
\begin{equation}
\label{LowerBound}
l^* = \arg\min_{l}\{\sum_{i=1}^l \rho(i) \ge \frac{1}{2} \} -1.
\end{equation}
In fact, $l^*$ denotes the largest integer that makes
\begin{equation}
\nonumber
\sum_{i\in \mathcal{S}} \rho(i)< \frac{1}{2},
\end{equation}
for an arbitrary index set $\mathcal{S}\subset [m]$ with $ |\mathcal{S}| \le l^*$.
If sequence $\{\rho(1),...,\rho(m)\}$ is not originally in a non-increasing order, one can sort it in a non-increasing order. And such a sort operation does not affect the conclusion in corollary 1.

According to corollary 1, we give a corresponding algorithm in table I to compute the exact value of $l^*$ when given a sensing matrix $A$.
The computational complexity of the presented algorithm is $O(m^2)$, which is the same as that of computing the coherence of $A$.
\begin{table}
  \centering                
\begin{tabular}{c}
\hline
\textbf {Table I}:\quad Algorithm for Lower Bound of Sparsity Level \quad \\
\hline
\end{tabular}
\begin{tabular}{l|l}
Input: & $\{\alpha_i\in R^{n}|\ i\in [m]\} $\\
\hline
\textit{Step 1}: &  \\
&$\textbf{for} \ {i = 1\ :\ m}$ \ \textbf{do}\\
& \quad compute $\{c_{ij} | j\in [m]\}$ with $c_{ij} = <\alpha_i, \alpha_j>$\\
&\quad $\nu(i)= \max_{j} \{\frac{|c_{ij}|}{c_{ii}} |\ j\neq i, \ j \in [m]\}$\\
&\quad $\rho(i) = \frac{\nu(i)}{1+ \nu(i)}$\\
&\quad $i = i+1$\\
&\textbf{end for}\\ 
\textit{Step 2}: &\\ 
&\textbf{Sort} $\{\rho(i) | \ i\in [m]\}$ in a non-increasing order\\
&\quad \quad  and get \{$\rho(i_1),...,\rho(i_{m})$ \}\\
\textit{Step 3}: & \\ 
& \textbf{Find} the smallest number $k$ that makes\\
&\quad \quad $\sum_{j = 1}^k \rho(i_j) \ge \frac{1}{2}$\\
&\textbf{Set} $l^* = k-1$\\ \hline
Output &  $l^*$\\ \hline
\end{tabular}      
\end{table}

\section{Discussions}
\subsection{Relation to Coherence}
In this subsection, we show that ASF improves on the bound of $SL$ derived by coherence. To make the comparison between ASF and coherence straightforward, we first assume $A$ consists of $l_2$ normalized columns. Therefore we have
\begin{equation}
\nonumber
\begin{split}
\nu(i)\le \mu, \ \forall i\in [m],
\end{split}
\end{equation}
which means
\begin{equation}
\nonumber
\rho(i)=\frac{\nu(i)}{1+\nu(i)}\le \frac{\mu}{1+\mu}.
\end{equation}
Then we get
\begin{equation}
\nonumber
\max_{|\mathcal{S} |=k}\rho(\mathcal{S})\le k\frac{\mu}{1+\mu}.
\end{equation}
According to \cite{gribonval2003sparse} and \cite{donoho2003optimally}, to ensure correct $l_1$ recovery of a sparse vector, the sparsity $k$ should be less than $\frac{1}{2}(1+\mu^{-1})$, which guarantees that
\begin{equation}
\nonumber
\max_{|\mathcal{S} |=k}\rho(\mathcal{S})<\frac{1}{2}.
\end{equation}
Above inequality implies that $l^*\ge k$. Therefore, we can say ASF is a finer measure for sparsity estimation than coherence.
It is important to note that, different from coherence which only considers the angle information between columns in a sensing matrix, ASF takes into account angle information as well as scale information of each column, which can improve ASF's performance in sparsity estimation.
\subsection{Relation to RIC}
RIP condition of order $s$ is commonly understood as a measure of "overall conditioning" of the set of $n\times s$ submatrices of $A$.
For simplicity, we first make several definitions:
Let $\mathcal{T}$ be an arbitrary set that has elements in $[m]$ with its cardinality equal to or less than $s$. Let $A_\mathcal{T}$ denote a $n\times s$ submatrix of $A$.
Then we can get two critical sets as below:
\begin{equation}
\nonumber
\begin{split}
\mathcal{T}_1 = \arg \max_{|\mathcal{T}|= s} \Lambda_{max}(A_{\mathcal{T}}^TA_{\mathcal{T}}) \\
\mathcal{T}_2 =\arg \min_{|\mathcal{T}|= s} \Lambda_{min}(A_{\mathcal{T}}^TA_{\mathcal{T}}).
\end{split}
\end{equation}
Then we set $k_{max}, k_{min}$ as
\begin{equation}
\label{ConditionNumbers2}
\begin{split}
k_{max} = \Lambda_{max}(A_{\mathcal{T}_1}^TA_{\mathcal{T}_1})\\
k_{min} = \Lambda_{min}(A_{\mathcal{T}_2}^TA_{\mathcal{T}_2}).
\end{split}
\end{equation}
Since the RIP condition does not hold the homogeneity property \cite{foucart2009sparsest}, to avoid this problem, we use the following condition instead:
\begin{equation}
\nonumber
k_{min}\|x\|_2^2 \le \|Ax\|_2^2 \le k_{max}\|x\|_2^2, \ \forall \ \|x\|_0\le s
\end{equation}
where $k_{max}>k_{min} >0$ and one can check parameter $\frac{k_{max}}{k_{min}}$ holds the homogeneity property. Obviously RIP condition is a special case of the above formulation by setting $k_{max} = 1+\sigma_s$ and $k_{min} = 1-\sigma_s$ where $\sigma_s \in (0,1)$.

\noindent\textbf{Lemma 2} Set $C = A^TA$. $k_{max}$ and $k_{min}$ are defined in ($\ref{ConditionNumbers2}$). There exist indices $h$ and $l$ such that
\begin{equation}
\nonumber
\begin{split}
c_{hh}- k_{min} \le \sum_{j\neq h, j\in \mathcal{T}_1} |c_{hj}|\le (s-1) |c_{hk}|\\
k_{max}- c_{ll} \le \sum_{j\neq l, i\in \mathcal{T}_2} |c_{lj}|\le (s-1) |c_{lf}|,
\end{split}
\end{equation}
where $|c_{hk}| = \max_{j\neq h}\{|c_{hj}|\}$ and $|c_{lf}| = \max_{j\neq l}\{|c_{lf}|\}$.
This Lemma can be easily derived using the Gersgorin Disc Theorem.

Based on Lemma 2, we get the following theorem:

\noindent \textbf{Theorem 2} With $k_{max}$ and $k_{min}$ defined in ($\ref{ConditionNumbers2}$), we have
\begin{equation}
\nonumber
\frac{k_{min}}{k_{max}} \ge \min_{h\neq l} \frac{c_{hh}}{c_{ll}} - 2(s-1)\max_i\nu(i)
\end{equation}

\noindent \textbf{Proof} According to lemma 2, we have
\begin{equation}
\nonumber
\begin{split}
&c_{hh}- k_{min}\le (s-1) |c_{hk}|\\
&\frac{c_{hh}}{c_{ll}}k_{max}- c_{hh} \le \ (s-1)\frac{c_{hh}}{c_{ll}} |c_{lf}|,
\end{split}
\end{equation}
then we get
\begin{equation}
\nonumber
\begin{split}
\frac{c_{hh}}{c_{ll}}-\frac{k_{min}}{k_{max}} &\le (s-1) \frac{|c_{hk}|}{k_{max}} + (s-1) \frac{|c_{lf}|}{c_{ll}} \frac{c_{hh}}{k_{max}}\\
&\le 2(s-1)\max_i \nu(i)
\end{split}
\end{equation}
with $\nu(i)$ is defined in $(\ref{Property1})$.
Furthermore, we get
\begin{equation}
\label{RIPConclusion}
\frac{k_{min}}{k_{max}} \ge \min_{i\neq j}\frac{c_{ii}}{c_{jj}} -2(s-1)\max_{i}\nu(i).
\end{equation}
$ \blacksquare$

By combining $k_{max}$, $k_{min}$ with the definition of RIP, we have:
\begin{equation}
\nonumber
\frac{k_{min}}{k_{max}} = \frac{1-\sigma_s}{1+\sigma_s}.
\end{equation}
According to $(\ref{RIPConclusion})$, we further have
\begin{equation}
\nonumber
\frac{1-\sigma_s}{1+\sigma_s} \ge \min_{i\neq j}\frac{c_{ii}}{c_{jj}} -2(s-1)\max_{i}\nu(i)
\end{equation}
Given $\max_{i}\nu(i) <\frac{1}{2(s-1)}\min_{i\neq j}\frac{c_{ii}}{c_{jj}}$, we have
\begin{equation}
\label{RIC estimation}
\nonumber
\sigma_s \le \frac{1-(\min_{i\neq j}\frac{c_{ii}}{c_{jj}} -2(s-1)\max_{i}\nu(i))}{1+(\min_{i\neq j}\frac{c_{ii}}{c_{jj}}-2(s-1)\max_{i}\nu(i))}.
\end{equation}
 When given a specific number $t\in (0,1)$, if we have
\begin{equation}
\nonumber
s <  \frac{1}{2\max_i \nu(i)}(\min_{i\neq j}\frac{c_{ii}}{c_{jj}}-\frac{1-t}{1+t}) +1,
\end{equation}
then it is guaranteed that $\sigma_s<t$. On the other hand, $t = \frac{1}{3}$ is a sharp bound to guarantee correct recovery of any $s$-sparse vector via $l_1$ minimization \cite{cai2013sharp}. Therefore we can get another lower bound for $SL$ as follows:
\begin{equation}
SL\ge \frac{1}{2\max_i \nu(i)}(\min_{i\neq j}\frac{c_{ii}}{c_{jj}}-\frac{1}{2}).
\end{equation}

\section{Conclusion}
In this paper, we successfully developed a new ASF for analyzing the sparsity recovery behavior of $l_1$ minimization. We prove that ASF provides a sharper bound for sparsity level of a sensing matrix than coherence. Also we further analyze the underlying relationship between ASF and RIC and derive an alternative RIC-based bound for sparsity level. ASF may find applications in areas like sparse coding and dictionary learning, coherent sampling and so on. ASF may provide hints on how to re-scale the columns of a sensing matrix in order to enhance sparsity. In future work, we are interested in applying ASF to analyze the stable recovery ability of the $l_1$ minimization problem in noisy environment.

\bibliographystyle{IEEEbib}
\bibliography{strings,refs}

\end{document}